\documentclass[12pt]{article}
\usepackage[cp1251]{inputenc}
\usepackage[T2A]{fontenc}
\usepackage[english]{babel}
\usepackage{mathrsfs}

\usepackage{mathtext,amsmath, amstext, amsgen, amsthm, mathrsfs, amsfonts, amssymb}
\newtheorem{thm}{Theorem}

\usepackage{setspace}
\begin{document}
\begin{center}

{\bf Global deformations of a Lie algebra of type $\bar{A_5}$}

{\bf  N.G. Chebochko, M.I. Kuznetsov}\\
{\it  National Research University Higher School of Economics, Nizhni Novgorod State University, Nizhni Novgorod, Russian Federation}

\end{center}
\textbf{\textit{It is shown that the orbits of the space of local deformations of the Lie algebra $\bar{A_5}$ over an algebraically closed field $K$ of characteristic 2 with respect to the automorphism group $\mathrm{PGL} (6)$ correspond to $\mathrm{GL} (V)$-orbits of tri-vectors of a 6-dimensional space. For local deformations corresponding to tri-vectors of rank $\rho <6$, integrability is proved and global deformations are constructed.}}

\textbf{\textit{Key Words:}} classical Lie algebra, field of characteristic 2, Lie algebra cohomology, integrable cocycle, deformations of Lie algebras.

\textbf{\textit{2010 Mathematics Subject Classification:}} Primary 17B50; 17B56; Secondary 17B66; 14M15.

\section{Introduction}

Deformations of Lie algebras are of considerable interest in the study of simple modular Lie algebras, especially in the case of fields of low characteristic. Classical Lie algebras over fields of characteristic zero and characteristic $p > 3$ are rigid (see \cite{rudakov}). In \cite{natimi}, \cite{kiril} it is proved that over fields of characteristic $p > 2$ only algebras of type $C_2$ have non-trivial deformations, their description is given in \cite{kostrikin},
  \cite{kostrikinkuz}. At present there does not seem to exist any complete description of deformations of classical Lie algebras of characteristic 2 and Lie algebras of their derivations. Spaces of local deformations of classical Lie algebras with a homogeneous root system are found in \cite{chebochko}. Among the Lie algebras of series $A$, only the algebras $\bar{A}_3, ~\bar{A}_5$ have nontrivial local deformations. Global deformations $\bar{A}_3$ are described in \cite{CK}. In \cite{CKK}, isomorphisms are constructed between deformations $\bar A_3$ and the known simple 14-dimensional Lie algebras of characteristic 2.

In this paper, we study the global deformations of a simple Lie algebra $L$ of type $\bar{A_5} = A_5 / Z$ over an algebraically closed field $K$ of characteristic 2. Here, as in \cite{kostrikinkuz} and \cite{CK}, the consideration is based on the study of the orbits of the space of local deformations $H^2(L, L)$ with respect to the automorphism group $G$ of the algebra $L$, since cocycles belonging to the same orbit define isomorphic families of deformations. In our case, $G = \mathrm{PGL} (6) \subset \mathrm{Aut} (L)$. In the paper an isomorphism of $\mathrm{GL}(6)$-module $H^2(L, L)$ and the module of tri-vectors shifted by Frobenius morphism and character $\det^{- 1}$ is constructed. This makes possible to reduce the classification of $G$-orbits of $H^2(L, L)$ to the classification of $\mathrm{GL} (V)$-orbits of tri-vectors in the 6-dimensional space $V$. In the case of the field of complex numbers, this problem was solved by W. Reichel in 1907  \cite{Reichel} (see \cite{Gur0}).
Namely, there are 5 orbits of tri-vectors with the following representatives:

(I) 0;

(II) $e_1\wedge e_2\wedge e_3$;

(III) $e_1 \wedge(e_2\wedge e_3 + e_4\wedge e_5)$;

(IV) $e_1\wedge e_5\wedge e_6 + e_2\wedge e_6\wedge e_4 + e_3\wedge e_4\wedge e_5$;

(V) $e_1\wedge e_2\wedge e_3 + e_4\wedge e_5\wedge e_6$.

It means that for a tri-vector $w$ there exists a basis of  $V$ with respect to which $w$ is written in one of the five canonical forms (I) – (V). Orbits (I) - (V) exist over fields of any characteristic. Orbit (V) is that of general position (see \cite{VinPopov}, The Basic Table). The set of non-trivial tri-vectors of rank $\rho <6 $ consists of orbits (II) and (III). Note that for $n $-dimensional space the classification is obtained by J.A. Schouten (n = 7, 1931, \cite{Schouten}, see \cite{Gur0}), G.B. Gurevich (n = 8, 1934-35, \cite{Gur} - \cite{Gur2}), E.B. Vinberg and A.G. Elashvili (n = 9, 1978, \cite{VinElas}).
For short, we define the rank of local deformation as the rank of the corresponding tri-vector. In the paper, the integrability of local deformations of rank $\rho <6 $ is proved and global deformations are constructed. As a result, we obtain simple Lie algebras of types $\bar{A_5}$ (II), $\bar{A_5}$ (III). The case $\rho = 6$ remains open.

Now we introduce the basic definitions. Let $L$ be a Lie algebra over a field $K, ~ F=  K((t)), ~L_F=F\otimes _KL$.  Lie algebra $L_F$ with multiplication
$f_t(x, y) = [x, y] + t\phi_1(x,y)+t^2\phi_2(x, y)+\ldots$, where $\phi_i$ are bilinear mappings over $K$, is called a global deformation of the Lie algebra $L$. In particular, $\phi_1$ is a cocycle from $Z^2(L, L)$.  Cocycles from one class of cohomology give equivalent deformations, therefore, to find ones, we consider  $H^2(L, L)$ as the space of local deformations.
Let $ R $ be a root system of type $ A_5 $ with basis
$ \{\alpha_1 = \varepsilon_1- \varepsilon_2, ~ \alpha_2 = \varepsilon_2- \varepsilon_3, ~ \alpha_3 = \varepsilon_3- \varepsilon_4, ~ \alpha_4 = \varepsilon_4-
\varepsilon_5, ~ \alpha_5 = \varepsilon_5- \varepsilon_6 \} $, $A$ be
Lie algebra of type $A_5$ over a field $K$ of characteristic 2, $ \{H _{\alpha_i} (i = 1, \ldots, 5), E_\alpha (\alpha \in R) \}$ be a Chevalley basis of $L$, $ \dim A = 35 $. The algebra $A$ has a one-dimensional center $Z$ with basis $ H_ {\alpha_1} + H_ {\alpha_3} + H_ {\alpha_5} $. We denote the quotient algebra $A/Z$ by $L$, and the basis of $L$ is denoted as in $A$. The algebra $A_5$ is seen as the Lie algebra $\mathrm{sl} (V), ~\dim V = 6 $, and the algebra $L$ as $\mathrm{sl} (V) / Z $. The adjoint representation gives the action of the group $\mathrm{GL} (V)$ by automorphisms of the Lie algebras $\mathrm{sl} (V)$ and $L = \mathrm{sl} (V) / Z$. In this paper, the roots and weights of cochains of the Lie algebra $L$ are considered with respect to a maximal torus of $\mathrm{GL} (V)$.

\section {Group $H ^ 2 (L, L)$ and  tri-vectors}

The second cohomology group $ H^2 (L, L) $ is found in \cite{chebochko}. It is proved
that $ \dim H ^ 2 (L, L) = 20 $, the weights of $ H ^ 2 (L, L) $ are conjugate with respect to the Weyl group with ${\alpha_1} + {\alpha_3} + {\alpha_5} $, any weight subspace $ H_ \mu ^ 2 (L, L) $ is one-dimensional and having a basis cocycle
$ \psi_ \mu = \sum \limits_{- \gamma- \delta + \mu \in R} E ^ * _ {- \gamma} \wedge
E ^ * _ {- \delta} \otimes E _ {- \gamma- \delta + \mu} $.
The root vector $ E _ {\alpha}, ~ \alpha = \varepsilon_i- \varepsilon_j $ of the Lie algebra $ L $ will be identified with  $ E _ {\varepsilon_i- \varepsilon_j} = e_i \otimes f_j \in \mathrm{sl} (V)$, where $ \{e_i \}$ is  a basis of $ V $, $ \{f_i \} $ is the dual basis of  $V^*$. In fact, we employ the sections of canonical
projections of quotient  modules by submodules over $\mathrm{GL}(V)$ and morphisms of restrictions on submodules. For example, for the natural section $s:L \longrightarrow \mathrm{sl}(V)$ we have $s(E_\alpha)=e_i\otimes f_j, ~\alpha =\varepsilon_i - \varepsilon_j, ~s(\bar{H}_i)=H_i, ~i=1, 2, 3, 4$. As a result, we obtain the embedding of $C^2(L, L)$ into $\mathrm{gl}(V)^*\wedge \mathrm{gl}(V)^*\otimes \mathrm{gl}(V)$.
The cocycle $\psi_\mu $ is easy to find, e.g., the weight $ \mu = {\alpha_1} + {\alpha_3} + {\alpha_5} = \varepsilon_1- \varepsilon_2 + \varepsilon_3- \varepsilon_4 +
\varepsilon_5- \varepsilon_6 $ may be put as the sum of the three roots in the following ways: $$ (\varepsilon_1- \varepsilon_2) + (\varepsilon_3- \varepsilon_4) +
(\varepsilon_5- \varepsilon_6) $$$$ (\varepsilon_1- \varepsilon_2) + (\varepsilon_3- \varepsilon_6) +
(\varepsilon_5- \varepsilon_4) $$$$ (\varepsilon_1- \varepsilon_4) + (\varepsilon_3- \varepsilon_2) +
(\varepsilon_5- \varepsilon_6) $$$$ (\varepsilon_1- \varepsilon_4) + (\varepsilon_3- \varepsilon_6) +
(\varepsilon_5- \varepsilon_2) $$$$ (\varepsilon_1- \varepsilon_6) + (\varepsilon_3- \varepsilon_2) +
(\varepsilon_5- \varepsilon_4) $$$$ (\varepsilon_1- \varepsilon_6) + (\varepsilon_3- \varepsilon_4) +
(\varepsilon_5- \varepsilon_2). $$

The remaining weights are arranged similarly as the sum of summands $ \varepsilon_1, \varepsilon_2, \varepsilon_3, \varepsilon_4,
\varepsilon_5, \varepsilon_6 $, where there are three terms with the $ + $ sign and three other with the $ - $ sign.

Obviously, the weight $ \pm \varepsilon_1 \pm \varepsilon_2 \pm \varepsilon_3 \pm \varepsilon_4 \pm
\varepsilon_5 \pm \varepsilon_6 $ cannot be represented as the sum of the two roots (the roots in $ A_5 $ are $ \varepsilon_i- \varepsilon_j $), therefore $ H_ \mu ^ 2 (L, L) = Z_ \mu ^ 2 (L , L) $.

Here, we give an implementation of $ H ^ 2 (L, L) $, a sketch being presented in \cite{ckperm}. The aim is to obtain the realization of $\mathrm{GL}(V)$-module $H^2(L, L)$ as a quotient of $\mathrm{GL}(V)$-modules isomorphic to a shifted module of tri-vectors of $V$.
Given this identification, we can assume that $ H ^ 2 (L, L) $ is a subspace in
$$ ((V \otimes V ^ *) ^ * \wedge (V \otimes V ^ *)^*) \otimes V \otimes V ^ *. $$
For any space $ U $, the tensor degree $ U ^ {\otimes k} $ can be considered
as the space of multilinear functions on the space $ U ^ * $. The application of the operation of alternating over the field of characteristic 2 can be written as
$$ A (v_1 \otimes \ldots \otimes v_k) = \sum \limits _ {\pi \in S_k} v _ {\pi (1)} \otimes \ldots \otimes v _ {\pi (k)} $$
(here we do not divide by $ k! $).
 In particular, for the alternating function, we have $ v_1 \wedge v_2 \wedge v_3 = A (v_1 \otimes v_2 \otimes v_3) $. Thus, we consider $ \Lambda ^ 3V $ to be embedded in $ V ^ {\otimes 3} $. Obviously, this embedding is equivariant with respect to
$ \mathrm{GL} (V) $. Similarly, $ ((V \otimes V ^ *) ^ * \wedge (V \otimes V ^ *)^*) \otimes V \otimes V ^ * $ is embedded in $ V ^ {\otimes 3} \otimes V ^ { * \otimes 3} $ as a space of functions skew symmetric with respect to the first two arguments from $ V ^ * $ and those from $ V $. Therefore, in $ V ^ {\otimes 3} \otimes V ^ {* \otimes 3} $ we have
$$ \Lambda ^ 3V \otimes \Lambda ^ 3V ^ * \subset ((V \otimes V ^ *) ^ * \wedge (V \otimes V ^ *)^*) \otimes V \otimes V ^ *. $$
For short, we do not use the $ \wedge $ sign for external multiplication.
Note that in the basis $ \{e_i \} $ of the space $ V $, the element $e_i e_j e_k \otimes f_s f_q f_r$ as the element $ V^{\otimes 3} \otimes V^{*\otimes 3} $ is equal to $ A (e_i \otimes e_j \otimes e_k) \otimes A (f_s \otimes f_q \otimes f_r) $, being the sum of 36 basic terms. Given the form of the base cocycles $ \psi_ \mu $ and the embeddings described above, we obtain the following realization of the space $H ^ 2 (L, L) $:

$$ H ^ 2 (L, L) = <e_i e_j e_k \otimes f_s f_q f_r, ~ {i, j, k} \cup {s, q, r} = {1, 2, \ldots, 6}> \subset \Lambda^3V \otimes \Lambda^3V^*. $$

 We have constructed a natural realization of $ H ^ 2 (L, L)$, which, however, is not invariant with respect to the action of  $ \mathrm{GL} (V) $ by automorphisms of the algebra $ L $, despite the invariance of embeddings of various tensor spaces. The reason is that we identified $ H ^ 2 (L, L) $ with the direct sum of some spaces $ Z ^ 2 _ {\mu} (L, L) $, while $ H ^ 2 = Z ^ 2 / B ^ 2 $.
The group $G = \mathrm{PSL} (V)$ is a connected component of the automorphism group of a Lie algebra $L$. We identified the space $ H ^ 2 (L, L) $ with the subspace in $G$-module $ \Lambda ^ 3V \otimes \Lambda ^ 3V ^ * $. The map of exterior multiplication $ \Lambda ^ 3V \times \Lambda ^ 3V \longrightarrow \Lambda ^ 6V $ is $ \mathrm{GL} (V) $-invariant pairing into the one-dimensional module $ K _ {\delta} $, where $ \delta (g) = \det ( g).$ Therefore, $ \mathrm{GL} (V) $- module $ \Lambda ^ 3V ^ * $ is isomorphic to the module
$ K _ {\delta ^ {- 1}} \otimes  \Lambda ^ 3 (V) $. Obviously, the element $ f_l f_m f_n $ corresponds to the element $ e_i e_j e_k $, such that $ \{i, j, k, l, m, n \} = \{1, 2, \ldots, 6 \}. $ Further, the $ \mathrm{GL} (V) $ - module $ \Lambda ^ 3V \otimes \Lambda ^ 3V ^ * $ is isomorphic to the $ \mathrm{GL} (V) $-module $ K _ {\delta ^ {- 1}} \otimes \Lambda ^ 3V \otimes \Lambda ^ 3V $. On these modules, the center $\mathrm{GL} (V)$ acts trivially. Therefore, we have an isomorphism of $ G$-modules which the subspace we identified with $ H ^ 2 (L, L) $ maps to the subspace
$$ <e_i e_j e_k \otimes e_i e_j e_k, ~ i <j <k> \subset  K _ {\delta ^ {- 1}} \otimes \Lambda ^ 3V \otimes \Lambda ^ 3V. $$
Denote $ \Lambda ^ 3V $ by $ W $. The space $ W \otimes W $ is that of all bilinear functions on $ W ^ * $. Denote by $ S ^ {(2)} (W) $ the module of symmetric functions. Obviously,
$ \Lambda ^ 2W \subset S ^ {(2)} (W) \subset W \otimes W $ and all embeddings considered  are consistent with the action of $ \mathrm{GL} (V) $. Besides,
$$ S ^ {(2)} (W) = \Lambda ^ 2W \oplus <e_i e_j e_k \otimes e_i e_j e_k, ~ i <j <k>. $$
Therefore, we have the isomorphism of $ \mathrm{GL} (V) $-modules
$$ H ^ 2 (L, L) \cong K _ {\delta ^ {- 1}} \otimes (S ^ {(2)} (W) / \Lambda ^ 2W). $$
The space $ S ^ {(2)}(W) $ is isomorphic as $\mathrm{GL} (V) $-module to the space $ \mathfrak {m} _2 $ of homogeneous elements of the second degree of the divided powers algebra $ \mathcal {O}(W) $ (see \cite{KostrikinShafarevich}), here $ \mathfrak {m} $ is the maximal ideal of $ \mathcal {O} (W) $, while $ \Lambda ^ 2W $ is isomorphic to $ \mathfrak {m} ^ 2_2 $, consisting of homogeneous second degree elements of $\mathfrak{m}^2$. Let $ \{w_i \} $ be a basis of $ W $. The classes of elements $ w_i ^{(2)}, ~ i = 1, \ldots, 20 $  form the basis of the space $ \mathfrak{m} _2 / \mathfrak{m}^ 2_2 $.
Thus, we have obtained the following theorem:

\begin{thm} Let $L = \mathrm{sl} (V) / Z $ be a Lie algebra of type $ \bar{A_5} $ over an algebraically closed field $K$ of characteristic 2, $ W = \Lambda ^ 3V $.

(i) The group $ H^2 (L, L) $ is isomorphic to the $\mathrm{GL}(V) $-module $K_{\delta ^{- 1}}
\otimes \mathfrak{m}_2 / \mathfrak{m}^ 2_2 $, where $ \mathfrak{m}_2 $ is the second term of the standard grading of the maximal ideal of the divided powers algebra  $\mathcal{O}(W)$.

(ii) The action of $\mathrm{GL}(V)$ on $\mathfrak{m}_2 / \mathfrak{m}^ 2_2$ is distinguished from  the standard action on the space of tri-vectors of the space $ V $ by Frobenius morphism.

(iii) The orbits of the projective space $ P(H^2(L, L)$ with respect to the group $\mathrm{GL}(V)$ are in one-to-one correspondence with the $\mathrm{GL} (V) $-orbits of the space $ P(\Lambda ^3V)$.
\end{thm}

      The transition to the projective space in statement (iii) is explained by the factor
$ K _{\delta^{- 1}} $ in (i). Obviously, cocycles that differ by a constant factor lead to isomorphic global deformations.
The classification of tri-vectors of a 6-dimensional space over $ \mathbb {C} $ is known owing to W. Reichel (see \cite{Gur}). Any tri-vector can be reduced to one of  5 the canonical types listed in the introduction. Corresponding orbits take place over any algebraically closed field. The stabilizer in $\mathrm{gl}(V)$ of a tri-vector of the form (V) is equal to $ \mathrm{sl} (3) + \mathrm{sl} (3) $, therefore the corresponding orbit has dimension 20 and, therefore, is the orbit of general position (see \cite{VinPopov}, The Basic Table). According to the theorem there is a natural one-to-one correspondence between the elements of $H^2(L, L)$ and tri-vectors. Recall that we define {\it the rank $\rho$} of a cocycle as the rank of corresponding tri-vector.

\section{Global deformations of $\bar{A_5}$ with the local part of rank $\rho < 6$}

A necessary condition for the prolongation of the cocycle $ \psi $ to global
deformation is the triviality of the cocycle $ \psi \cup \psi $ in
$ H^3(L, L) $, where $ \psi \cup
\chi (x, y, z) = \psi (\chi (x, y), z) + \psi (\chi (y, z), x) + \psi (\chi (z, x), y ) $.
We also define $ [\psi, \chi] = \psi \cup \chi + \chi \cup \psi $.

From the definition of the natural action of $\mathrm{Aut} (L) $ on $Z^k (L, L)$ right away
follows that
for any $ g\in G $, $ g [\psi, \chi] = [g\psi, g\chi]$.

We investigate the prolongation of representatives of orbits of types (II) and (III).

The cocycle corresponding to tri-vector of the type (II) is as follows $ \psi = \psi _ {\varepsilon_1 + \varepsilon_2 + \varepsilon_3- \varepsilon_4- \varepsilon_5- \varepsilon_6}.$ The cocycle $ \psi $ has a nonzero value on pairs of basis vectors of the form $ E _ {- \varepsilon_i + \varepsilon_j} (i = 1,2,3, j = 4,5,6)$ only, and the image of the pair is some $ E_ { \varepsilon_i- \varepsilon_j} (i = 1,2,3, j = 4,5,6) $. Therefore, $ \psi (\psi (x, y), z) = 0 $, $ \psi $ extends to global deformation and $ [~, ~] + t \psi $ defines the Lie multiplication. The corresponding Lie algebras are denoted by $ \bar{A_5}(\mathrm{II}) $.

The cocycles corresponding to tri-vector of the type (III) is as follows $ \psi = \psi _ {\varepsilon_1 + \varepsilon_2 + \varepsilon_3- \varepsilon_4- \varepsilon_5- \varepsilon_6} + \psi _ {\varepsilon_1-\varepsilon_2-
\varepsilon_3 + \varepsilon_4 + \varepsilon_5- \varepsilon_6}$.  Let $ \psi_1 = \psi _ {\varepsilon_1 + \varepsilon_2 + \varepsilon_3- \varepsilon_4- \varepsilon_5- \varepsilon_6} $, $ \psi_2 = \psi _ {\varepsilon_1-\varepsilon_2-
\varepsilon_3 + \varepsilon_4 + \varepsilon_5- \varepsilon_6}$. Then $ \psi \cup
\psi = \psi_1 \cup
\psi_1 + [\psi_1, \psi_2] + \psi_2 \cup
\psi_2 = [\psi_1, \psi_2] $ (for all basic cocycles $ \psi \cup \psi = 0 $).

We introduce an abbreviated notation for 2 cochains: if $ \varphi (E _ {\varepsilon_i- \varepsilon_j}, E _ {\varepsilon_k- \varepsilon_s}) = E _ {\varepsilon_t- \varepsilon_r} $, then we think that the set $ [ \{\varepsilon_i- \varepsilon_j, \varepsilon_k- \varepsilon_s \}, \varepsilon_t- \varepsilon_r] $ enters the cochain $ \varphi $ and conditionally write the cochain $ \varphi $ as the sum of the sets included. In a similar way, we write 3-cochains.

The cocycle $ [\psi_1, \psi_2] $ has the weight $ 2 (\varepsilon_1- \varepsilon_6) = (\varepsilon_1 + \varepsilon_2 + \varepsilon_3- \varepsilon_4- \varepsilon_5- \varepsilon_6)+
(\varepsilon_1- \varepsilon_2-
\varepsilon_3 + \varepsilon_4 + \varepsilon_5- \varepsilon_6) $. Determine how many sets are included in $ \psi_1 \cup \psi_2 $. In order for the set to be in $ \psi_1 \cup \psi_2 $, the last element of the set included in $ \psi_2 $ must coincide with one of the first elements in the set of $ \psi_1 $. It is only possible for $ - \varepsilon_2 + \varepsilon_4 $, $ - \varepsilon_2 + \varepsilon_5 $, $ - \varepsilon_3 + \varepsilon_4 $, $ - \varepsilon_3 + \varepsilon_5 $. Consider, e.g., $ - \varepsilon_2 + \varepsilon_4 $. $\psi_2 $ includes $ [\{- \varepsilon_1 + \varepsilon_3, - \varepsilon_5 + \varepsilon_6 \}, - \varepsilon_2 + \varepsilon_4] $ and $ [\{- \varepsilon_1 + \varepsilon_6, \varepsilon_3-\varepsilon_5\},- \varepsilon_2 + \varepsilon_4]$. For the first set, there are 4 variants of the third element in the corresponding set in $\psi_1 \cup \psi_2 $: $-\varepsilon_1+\varepsilon_5, -\varepsilon_1+\varepsilon_6, -\varepsilon_3+\varepsilon_5, -\varepsilon_3+\varepsilon_6$. As for the second, there are 3 only, since there is a coincidence with $ - \varepsilon_1 + \varepsilon_6 $. There are only 7 options available. The situation is similar for the rest of sets. We get that $ \psi_1 \cup \psi_2 $ includes $ 7 \cdot4 = 28 $ sets totally. The same number of sets is included in $ \psi_2 \cup \psi_1 $. But $ \psi_1 \cup \psi_2 $ and $ \psi_2 \cup \psi_1 $ have matches. In each group of 7 elements described above there is an element of the form $ [\{- \varepsilon_1 + \varepsilon_6, \varepsilon_3- \varepsilon_5, - \varepsilon_3 + \varepsilon_5 \}, \varepsilon_1- \varepsilon_6] $. It is repeated in $\psi_2\cup \psi_1 $. Therefore, $[\psi_1, \psi_2] $  includes 48 sets.

Consider the cochain \\$\varphi=[\{-\varepsilon_1+\varepsilon_5, -\varepsilon_5+\varepsilon_6\}, \varepsilon_1-\varepsilon_6]+[\{-\varepsilon_1+\varepsilon_5, -\varepsilon_1+\varepsilon_6\}, \varepsilon_5-\varepsilon_6]+\\+[\{-\varepsilon_1+\varepsilon_6, -\varepsilon_5+\varepsilon_6\}, \varepsilon_1-\varepsilon_5]+[\{-\varepsilon_4+\varepsilon_6, -\varepsilon_1+\varepsilon_4\}, \varepsilon_1-\varepsilon_6]+\\+[\{-\varepsilon_1+\varepsilon_4, -\varepsilon_1+\varepsilon_6\}, \varepsilon_4-\varepsilon_6]+[\{-\varepsilon_1+\varepsilon_6, -\varepsilon_4+\varepsilon_6\}, \varepsilon_1-\varepsilon_4]$. We can show that $d\varphi$ also includes 48 sets. Note that $d\varphi(E_{-\varepsilon_1+\varepsilon_5}, E_{-\varepsilon_5+\varepsilon_6}, E_{-\varepsilon_1+\varepsilon_6})=0$ and $d\varphi(E_{-\varepsilon_1+\varepsilon_4}, E_{-\varepsilon_4+\varepsilon_6}, E_{-\varepsilon_1+\varepsilon_6})=0$, so when counting the sets included in the differential of the parts of $ \varphi $ we exclude them. As a result we obtain $ H _ {\alpha_1} + H _ {\alpha_3} + H _ {\alpha_5} $ which is the central element of $ A_5 $.
The number of remaining sets in $d[\{-\varepsilon_1+\varepsilon_5, -\varepsilon_5+\varepsilon_6\}, \varepsilon_1-\varepsilon_6]$ is 14. In fact, there are the sets $[\{-\varepsilon_1+\varepsilon_5, -\varepsilon_5+\varepsilon_6, -\varepsilon_1+\varepsilon_i\}, \varepsilon_i-\varepsilon_6]$ ($i=2, 3, 4$) the sets $[\{-\varepsilon_1+\varepsilon_5, -\varepsilon_5+\varepsilon_6, \varepsilon_6-\varepsilon_i\}, \varepsilon_1-\varepsilon_i]$ ($i=2, 3, 4$) the sets $[\{-\varepsilon_1+\varepsilon_i, -\varepsilon_i+\varepsilon_5, -\varepsilon_5+\varepsilon_6\}, \varepsilon_1-\varepsilon_6]$ ($i= 2, 3, 4, 6$) and the sets $[\{-\varepsilon_1+\varepsilon_5, -\varepsilon_5+\varepsilon_i,-\varepsilon_i+\varepsilon_6\}, \varepsilon_1-\varepsilon_6]$ ($i=1, 2, 3, 4$), all in all 14. The same number of sets can be observed in $d[\{-\varepsilon_4+\varepsilon_6, -\varepsilon_1+\varepsilon_4\}, \varepsilon_1-\varepsilon_6]$.
 In the remaining parts of $ d \varphi $ there are 13 sets, the total number being 80. In the differentials 1 and 4 as well as 3 and 6 parts of $ \varphi $ there are two matches, while in 2 and 5 parts being no matches. The remaining pairs have one match each. Thus, from 80 sets under consideration, in fact, 48 ones should be left (80-32 = 48).
 Further, it is easy to verify that if a set is included in $ [\psi_1, \psi_2] $, then it enters $ d \varphi $ as well. Therefore, $ [\psi_1, \psi_2] = d \varphi $.

 Since $ \varphi $ has the weight $ 2 (\varepsilon_1- \varepsilon_6) $, then $ \varphi \cup \psi_1 $ has the weight $ 3 \varepsilon_1 + \varepsilon_2 + \varepsilon_3- \varepsilon_4- \varepsilon_5-3\varepsilon_6$ which cannot de obtained as a sum of four roots (the sum of positive coefficients is 5). The cochain $\varphi \cup \psi_2$ has the weight $3 \varepsilon_1-\varepsilon_2- \varepsilon_3 + \varepsilon_4 + \varepsilon_5-3\varepsilon_6 $, which cannot be obtained as the sum of four roots, either. Finally, $\varphi \cup \varphi$ has the weight $ 4 (\varepsilon_1- \varepsilon_6) $ and cannot be obtained as the sum of four roots, if at least three of which are different. Therefore, $ \psi_1 + \psi_2 $ continues to global deformation and $ [~, ~] + t (\psi_1 + \psi_2) + t ^ 2 \varphi $ defines the Lie multiplication. The corresponding Lie algebras are denoted by $ \bar{A_5}(\mathrm{III}) $.
Thus, cocycles of the rank $\rho < 6$ continue to the global deformations $\bar{A_5}(\mathrm{II}) $ and $ \bar{A_5}(\mathrm{III}) $.

{\bf Acknowledgments}
The authors are very grateful to A.V. Balandin for useful discussions. The investigation is funded by RFBR according to research project N 18-01-00900.

\end{document}